\documentstyle[12pt]{amsart}
\begin{document}
\baselineskip=16.5pt

\title[Nakano positivity and the
$L^2$-metric]{Nakano positivity and the $L^2$-metric
on the direct image of an adjoint positive line bundle}

\author[I. Biswas]{Indranil Biswas}

\address{School of Mathematics, Tata Institute of Fundamental
Research, Homi Bhabha Road, Bombay 400005, INDIA}

\email{indranil@@math.tifr.res.in}

\author[Ch. Mourougane]{Christophe Mourougane}

\address{Institut de Math\'ematiques de Jussieu,
Bo\^{\i}te 247, 4, place Jussieu, 75 252 Paris Cedex 05, FRANCE}

\email{mourouga@@math.jussieu.fr}

\date{}

\maketitle

\section{Introduction}

A holomorphic vector bundle
$E$ over a complex projective manifold $X$ is called {\it ample}
if the tautological line bundle ${\cal O}_{{\Bbb P}(E)}(1)$
over the projective bundle ${\Bbb P}(E)$ of hyperplanes
in $E$ is ample. This notion of ampleness was introduced by
R. Hartshorne in \cite{Ha}.

On the other hand, P.A. Griffiths in \cite{Gr}
introduced an analytic notion of positivity of a vector bundle.
A holomorphic vector bundle $E$ over $X$ is called
{\it Griffiths positive} if it admits a
Hermitian metric $h$ such that the curvature
$C_h(E) \in C^{\infty}(X,\, {\Omega}^{1,1}_X\bigotimes
\mbox{End}(E))$ of the Chern connection on $E$ for the
Hermitian structure $h$ has the property that for every
$x\in X$ and every
nonzero holomorphic tangent vector $0\neq v \in
T_x X$, the endomorphism $\sqrt{-1}\cdot C_h(E)(x)(v, {\bar v})$
of the fiber $E_x$ is positive definite with respect to $h$.

If $h$ is a Griffiths positive Hermitian metric on $E$,
then the Hermitian metric on ${\cal O}_{{\Bbb P}(E)}(1)$
induced by $h$ has positive curvature. Therefore, $E$ is ample
by a theorem due to Kodaira.
An ample line bundle admits a Griffiths positive Hermitian
metric. Also, an ample vector bundle on a Riemann
surface is Griffiths positive, which was proved by H. Umemura
in \cite{Um}. However,
the question posed by Griffiths \cite[page 186, (0.9)]{Gr} asking
whether every ample vector bundle is Griffiths positive is yet
to be settled.

The notion of Griffiths positivity was strengthened by S. Nakano. A
holomorphic Hermitian vector bundle $(E,h)$ is called {\it Nakano
positive} if $\sqrt{-1}\cdot C_h(E)$ is a positive form on $TX
\bigotimes E$.
A Nakano positive vector bundle is clearly Griffiths positive.
In the other direction, H. Skoda and J.-P. Demailly proved that
if $E$ is Griffiths positive, then $E\otimes \det E$ is Nakano
positive \cite{DS}, $\det E$ being the line bundle given by the top
exterior power of $E$.

Our aim here is to establish the following analytic property
of an ample vector bundle (Theorem 3.1).

\medskip
\noindent {\bf Theorem\, A.}\, {\it Let $E$ be an ample vector
bundle of rank $r$ over a projective manifold $X$. The vector bundle
$S^k(E)\bigotimes \det E$ is Nakano positive for every
$k\geq 0$, where $S^k(E)$ denotes the $k$-th symmetric
power of $E$. More generally, for every decreasing sequence
$\lambda \in {\Bbb N}^r$ of height (= number
of positive elements) $l$, the vector bundle
${\Gamma}^{\lambda} E\bigotimes (\det E)^{\otimes l}$ is Nakano
positive; ${\Gamma}^{\lambda} E$ is the vector bundle corresponding
to the irreducible representation of $GL(r, {\Bbb C})$
defined by the weight $\lambda$. So, in particular,
${\Gamma}^{\lambda} E\bigotimes (\det E)^{\otimes l}$ is Griffiths
positive.}
\medskip

In the special case where $X$ is a toric variety or an abelian
variety, the Griffiths positivity of such vector bundles associated
to an ample bundle was proved in \cite{Mo}.

The above Theorem A is obtained as consequence of a result on the
Nakano positivity of the direct image of an adjoint positive
line bundle. This result on Nakano positivity of direct image
will be described next.

Let $L$ be a holomorphic line bundle, equipped with a Hermitian
metric $h$, over a connected projective manifold $M$.
Given any section $t \in H^0(M,\, L\bigotimes K_M)$,
where $K_M$ is the canonical bundle,
its conjugate $\overline{t}$ is realized as a section
of $L^*\bigotimes \overline{K_M}$ using $h$. Now, given another
section $s \in H^0(M,\, L\bigotimes K_M)$, consider
the top form on $M$
obtained from $s\wedge \overline{t}$ by contracting $L$ with $L^*$.
The $L^2$ inner product on $H^0(M, \, L\bigotimes K_M)$ is defined
by taking the integral of this form over $M$.

We prove the following theorem on the Nakano positivity of the
$L^2$ metric on a direct image (Theorem 2.3).

\medskip
\noindent {\bf Theorem\, B.}\, {\it Let $\psi \, : \, Y \,
\longrightarrow\, X$ be a holomorphically locally trivial
fiber bundle, where $X$ and $Y$ are connected projective manifolds,
and $H^1({\psi}^{-1}(x),\, {\Bbb Q}) \,=\,0$ for
some, hence every, point $x\in X$.
Let $L$ be a holomorphic line bundle over $Y$ equipped
with a positive Hermitian metric $h$. Then the
$L^2$-metric, defined using $h$, on the vector bundle
${\psi}_*(K_{Y/X}\bigotimes L)$
over $X$ is Nakano positive; $K_{Y/X}$ is the relative canonical
bundle.}
\medskip

We note that the condition of ampleness of $L$
in Theorem B ensures that
the direct image ${\psi}_*(K_{Y/X}\bigotimes L)$ is
locally free with the fiber of the corresponding vector bundle
over any point
$x\in X$ being $H^0({\psi}^{-1}(x),\,
(K_{Y/X}\bigotimes L)\big\vert_{{\psi}^{-1}(x)})$.

Theorem A is an immediate consequence of Theorem B applied to a
natural line bundle over the flag bundle $\psi : \,M_{\lambda}(E) \,
\longrightarrow \, X$ associated to an ample vector bundle $E$ by
a weight $\lambda$. In particular, setting
$Y = {\Bbb P}(E)$ and $L = {\cal O}_{{\Bbb P}(E)}(k+r)$ in Theorem B,
where $E$ is an ample vector bundle of rank $r$, the Nakano
positivity of $S^k(E)\bigotimes \det E$ is obtained. This
particular case corresponds to the weight $(k,0,0, \cdots ,0)$.

\medskip
\noindent {\it Acknowledgments:}\, Remarks 3.3 and 3.5 were
communicated to us by Jean-Pierre Demailly. We are very grateful
to him for this. We thank M.S. Narasimhan
for some useful discussions. We would like to
thank the International Centre for Theoretical
Physics, Trieste, for its hospitality. The first-named author
also thanks Joseph Le Potier and the Institut de Math\'ematiques
de Jussieu for hospitality. The second-named author thanks
T. Ohsawa and K. Yoshikawa for his visit to Nagoya University.

\section{Curvature of the $L^2$-metric}

We first recall the definition of Nakano positivity. Let $E$ be a
holomorphic vector bundle over $X$ equipped with a Hermitian
metric $h$. Let $C_h(E)$ denote
the curvature of the corresponding Chern connection on
$E$. Let ${\Theta}_h$ denote the unique sesquilinear form
on $TX\bigotimes E$ such that for any $v_1, v_2\in T_x X$ and
$e_1,e_2 \in E_x$, the equality
$$
{\Theta}_h (v_1\otimes e_1, v_2\otimes e_2) \, =\, \sqrt{-1}\cdot
\langle C_h(E)(v_1, \overline{v_2})e_1 \, , e_2\rangle_h
$$
is valid. The metric $h$ is called {\it Nakano positive} if 
the sesquilinear form ${\Theta}_h$ on $TX\bigotimes E$ is Hermitian,
i.e., it is positive definite.

Let $X$ and $Y$ be two connected complex
projective manifolds of dimension $m$ and $n$
respectively, and
let
$$
\psi \, : \, Y \, \longrightarrow \, X \leqno{(2.1)}
$$
be a holomorphic
submersion defining a holomorphically locally trivial
fiber bundle over $X$ of relative
dimension $f =n-m$. So, every point of $X$ has an analytic
open neighborhood $U$ such that ${\psi}^{-1}(U)$ is holomorphically
isomorphic to the trivial fiber bundle
$U\times F$ over $U$,
where $F$ is the typical fiber of $\psi$. The relative
canonical line bundle $K_Y\bigotimes {\psi}^*K^{-1}_X$
over $Y$ will be denoted by $K_{Y/X}$.

For any $x\in X$, the submanifold ${\psi}^{-1}(x)$ of $Y$ will
be denoted by $Y_x$.

Let $L$ be an ample line bundle over $Y$. The direct
image ${\psi}_*(K_{Y/X}\bigotimes L)$ is locally free
on $X$. Indeed, from the Kodaira vanishing theorem it
follows that all the higher direct images of $K_{Y/X}\bigotimes L$
vanish.
Let $V$ denote the vector bundle over $X$ given by this direct
image ${\psi}_*(K_{Y/X}\bigotimes L)$.

Fix a positive Hermitian metric $h$ on $L$.
For any point $x\in X$,
using the natural conjugate linear isomorphism of
${\Omega}^{f,0}_{Y_x}$ with ${\Omega}^{0,f}_{Y_x}$ and the
Hermitian metric $h$ on $L\vert_{Y_x}$, a conjugate linear
isomorphism between ${\Omega}^{f,0}_{Y_x}\bigotimes L\vert_{Y_x}$
and ${\Omega}^{0,f}_{Y_x}\bigotimes L^*\vert_{Y_x}$ is obtained.
We denote this conjugate linear isomorphism by $\iota$.

The $L^2$ Hermitian metric on the vector bundle
$V:= {\psi}_*(K_{Y/X}\bigotimes L)$ is defined by sending
any pair of sections
$$
t_1,\, t_2 \, \in \, V_x \, := \, H^0(Y_x\, , \, (K_{Y/X}\otimes
L)\big\vert_{Y_x})
$$
to the integral over $Y_x$ of the contraction of
$(\sqrt{-1})^{f^2}t_1$ with $\iota (t_2)$. In other words,
if $t_i = s_i \otimes d{\zeta}_1\wedge
d{\zeta}_2\wedge \cdots \wedge d{\zeta}_f$, where $\{{\zeta}_1,
\cdots , {\zeta}_f\}$ is a local holomorphic coordinate chart on the
fiber $Y_x$ and $s_i$, $i=1,2$, is a local section of $L$,
then the pairing
$$
\langle t_1,t_2\rangle \, := \, (\sqrt{-1})^{f^2}\int_{Y_x}
\langle s_1, s_2 \rangle_h \cdot d{\zeta}_1\wedge d{\zeta}_2\wedge
\cdots \wedge d{\zeta}_f\wedge d\overline{{\zeta}_1}\wedge
d\overline{{\zeta}_2} \wedge \cdots \wedge d\overline{{\zeta}_f}
\leqno{(2.2)}
$$
is the $L^2$ inner product of the two vectors $t_1$ and
$t_2$ of the fiber $V_x$. The top form
$$
\langle s_1, s_2 \rangle_h \cdot d{\zeta}_1\wedge d{\zeta}_2\wedge
\cdots \wedge d{\zeta}_f\wedge d\overline{{\zeta}_1}\wedge
d\overline{{\zeta}_2} \wedge \cdots \wedge d\overline{{\zeta}_f}
$$
clearly depends only on $t_1$ and $t_2$ and, in particular, it
does not depend on the choice
of the coordinate function $\zeta$; in other words, it is a
globally defined top form on $Y_x$. We note that the definition of
$L^2$-metric does not require a metric on $Y_x$.

Our aim is to compute the curvature of the
Chern connection on $V$ for the $L^2$ metric.

\medskip
\noindent {\bf Theorem\, 2.3.}\, {\it Assume that
$H^1(Y_x,\, {\Bbb Q}) \, =\, 0$ for
some, hence every, point $x\in X$.
Then the $L^2$-metric on the direct image
$V$ is Nakano positive.}
\medskip

{\it Proof.}\, Take a point $\xi \in X$. Let $\{z_1, z_2,
\cdots ,z_m\}$ be a holomorphic coordinate chart on $X$
around $\xi$ such that $\xi = 0$.

Let $r$ be the rank of the direct image $V$. Fix a normal frame
$\{t_1, t_2, \cdots ,t_r\}$ of $V$ around $\xi$
with respect to the $L^2$ metric. In other words,
$\{t_1, t_2, \cdots ,t_r\}$ is a holomorphic frame
of $V$ around $\xi$ such that for the function
$\langle t_{\alpha} , t_{\beta} \rangle$ around $\xi$ we have
\begin{eqnarray*}
\langle t_{\alpha} , t_{\beta} \rangle\big\vert_{z=0} & = &
{\delta}_{\alpha\beta}\\
\frac{\partial\langle t_\alpha , t_\beta\rangle}
{\partial z_i}\big\vert_{z=0} & = & 0
\end{eqnarray*}
for all $\alpha ,\beta \in [1,r]$ and all $i\in [1,m]$.
We note that the second
condition is equivalent to the condition that
$d \langle t_\alpha , t_\beta \rangle (0) =0$, where $d$ is the
exterior derivation.

Let ${\nabla}^{L^2}$ denote the Chern connection
on $V$ for the $L^2$ metric. Its curvature, which is
a ${\rm End}(V)$-valued
$(1,1)$-form on $X$, will be denoted by $C_{L^2}$. Take
any vector
$$
v\, = \, \sum_{i=1}^m\sum_{\alpha = 1}^r v_{i,\alpha}
\frac{\partial}{\partial z_i}\otimes t_\alpha \, \in \,
T_\xi X\otimes V_\xi \, .
$$
We wish to show that
$$
\sum_{i=1}^m \sum_{j=1}^m \sum_{\alpha = 1}^r \sum_{\beta = 1}^r
v_{i,\alpha} \overline{v_{j,\beta}}\sqrt{-1}\cdot \big\langle
C_{L^2}\left(\frac{\partial}{\partial z_i}, 
\frac{\partial}{\partial\overline{z_j}}\right)t_\alpha \, ,
t_\beta \big\rangle_{L^2} \, > \, 0 \leqno{(2.4)}
$$
with the assumption that $v \neq 0$. Now, we have
$$
\big\langle \sqrt{-1}\cdot C_{L^2}\left(\frac{\partial}{\partial z_i},
\frac{\partial}{\partial\overline{z_j}}\right)t_\alpha,
t_\beta\big\rangle_{L^2} (0)
\, =\, \frac{1}{\sqrt{-1}}\cdot\frac{{\partial}^2 \langle t_\alpha
\, , t_\beta\rangle}{{\partial z_i} \partial\overline{z_j}}(0) \, .
$$

For any $\alpha \in [1,r]$,
let $\widehat{t}_\alpha$ be the {\it unique} section of
$K_{Y/X}\bigotimes L$, defined
on an analytic open
neighborhood of $Y_{\xi} := {\psi}^{-1}(\xi)$, which is determined
by the condition that for any $x \in X$ in a sufficiently small
neighborhood of $\xi$, the restriction of $\widehat{t}_\alpha$ to
the fiber ${\psi}^{-1}(x)$ represents $t_\alpha (x)$,
the evaluation of the section $t_\alpha$ at
$x$. The holomorphicity of the section $t_\alpha$ of $V$ implies
that the section $\widehat{t}_\alpha$ of
$K_{Y/X}\bigotimes L$ is also holomorphic.

Fix a trivialization of the fiber bundle over a neighborhood
of $\xi$ defined by the projection $\psi$.
In other words, we fix an isomorphism of the
fiber bundle ${\psi}^{-1}(U) \,\longrightarrow \, U$ over some open
subset $U\subset X$ containing $\xi$ with the trivial fiber bundle
$U\times Y_{\xi}$ over $U$. Using this isomorphism, the relative
tangent bundle for the projection $\psi$
gets identified with a subbundle of the restriction of the tangent
bundle $TY$ to ${\psi}^{-1}(U)$.
Consequently, any relative differential form
on ${\psi}^{-1}(U)$ becomes a differential form on ${\psi}^{-1}(U)$.

If ${\theta}$ (respectively, ${\omega}$)
is a $L$-valued $(i_1,i_2)$-form
(respectively, $L$-valued $(j_1,j_2)$-form), i.e., a smooth
section of ${\Omega}^{i_1,i_2}_Y\bigotimes L$ (respectively,
${\Omega}^{j_1,j_2}_Y\bigotimes L$), then define $\{\theta ,\omega\}$
to be the $(i_1+j_2, i_2+j_1)$-form on $Y$ obtained from the section
$$
{\theta}\wedge \overline{\omega} \, \in \,
C^{\infty}({\Omega}^{i_1+j_2,i_2+j_1}_Y \otimes L\otimes L^*)\, ,
$$
where the conjugate of
$L$ has been identified with $L^*$ using the Hermitian structure $h$,
and then contracting $L$ with $L^*$.

In the above notation we have
$$
{\partial}_X \langle t_\alpha , t_\beta\rangle \, = \,
{\psi}_*({\partial}_Y\{\widehat{t}_\alpha ,
\widehat{t}_\beta\})\, , \leqno{(2.5)}
$$
where ${\psi}_*$ is the integration of forms along the
fiber (the Gysin map) and $\alpha , \beta\in [1,r]$.

It is easy to see directly that the right-hand side
of (2.5) does not depend on the
choice of the trivialization of the fibration defined by $\psi$.
Indeed, any two choices of the local
trivialization defines a homomorphism, over $U$, from the tangent
bundle $TU$ to the trivial vector
bundle over $U$ with the space of vertical vector fields $H^0(Y_{\xi},
\, TY_{\xi})$ as the fiber. This homomorphism is constructed by
taking the difference of the two horizontal lifts, given by the two
trivializations, of vector fields.
If the one-form in the right-hand side of
(2.5) for a different choice of trivialization is denoted by $\eta$,
then the one-form defined by the difference
$$
\eta \,- \,{\psi}_*(\partial_Y \{\widehat{t}_\alpha ,
\widehat{t}_\beta \})
$$
sends any tangent vector $v \in T_xX $ to
$$
\int_{{\psi}^{-1}(x)} L_{\hat{v}}\{\widehat{t}_\alpha ,
\widehat{t}_\beta\} \, ,
$$
where $ L_{\hat{v}}$ is the Lie derivative with respect to
the vertical vector field corresponding to $v$ for the given pair
trivializations of the fiber bundle. Now the identity
$L_{\hat{v}} = d\circ i_{\hat{v}} + i_{\hat{v}}\circ
d$ and the Stokes' theorem together
ensure that the right-hand side of (2.5) is independent of the
choice of trivialization of the fibration defined by $\psi$.

Let $\nabla$ denote the Chern connection of the holomorphic line
bundle $L$ equipped with the Hermitian metric $h$. Its curvature will
be denoted by $C_h$.

Since $\widehat{t}_\alpha$ is a holomorphic section for
every $\alpha \in [1,r]$, the equality
$$
{\psi}_*({\partial}_Y\{\widehat{t}_\alpha , \widehat{t}_\beta\}) \,
=\, {\psi}_*(\{\nabla\widehat{t}_\alpha , \widehat{t}_\beta \})
$$
is valid for all $\alpha ,\beta\in [1,r]$.

For any $i\in [1,m]$,
let $\frac{\widetilde\partial}{\partial z_i}$ denote vector field
on ${\psi}^{-1}(U)$ given by the lift of the local vector field
$\frac{\partial}{\partial z_i}$ on $U\subset X$
using the chosen trivialization of the fiber bundle.
The contraction of a one-form $\theta$ with a vector
field $\nu$ will be denoted by $(\theta ,\nu)$.

The above equality combined with (2.5) gives
$$
\frac{\partial}{\partial z_i} \langle t_\alpha , t_\beta\rangle \,=\,
\left({\psi}_*(\{\nabla_{\frac{\widetilde\partial}{\partial z_i}}
\widehat{t}_\alpha ,\widehat{t}_\beta\})\right)\big\vert_{Y_{\xi}}
$$
as any $\widehat{t}_\alpha$ is a {\it relative top form}
with values in $L$.
Now taking taking $\overline{\partial}_Y$ of this equality yields
$$
\frac{{\partial}^2}{\partial\overline{z_j}\partial z_i}
\langle t_\alpha , t_\beta \rangle \,=\,
\left(\overline{\partial}_X
{\psi}_*(\{\nabla_{\frac{\widetilde\partial}{\partial z_i}}
\widehat{t}_\alpha , \widehat{t}_\beta\}),
\frac{\partial}{\partial\overline{z}_j}\right) \,=\,
\left({\psi}_*(\overline{\partial}_Y
\{\nabla_{\frac{\widetilde\partial}{\partial z_i}}
\widehat{t}_\alpha , \widehat{t}_\beta\}),
\frac{\partial}{\partial\overline{z_j}}\right)\, .
$$
Since any $\widehat{t}_\alpha$
is holomorphic, the last term coincides with
$$
\left({\psi}_*(\{\nabla^{0,1}\nabla_{
\frac{\widetilde\partial}{\partial z_i}}
\widehat{t}_\alpha , \widehat{t}_\beta\}),
\frac{\partial}{\partial\overline{z_j}}\right)
\, +\, (-1)^f\cdot
\left({\psi}_*(\{\nabla_{\frac{\widetilde\partial}{\partial z_i}}
\widehat{t}_\alpha , \nabla^{1,0}\widehat{t}_\beta\}),
\frac{\partial}{\partial\overline{z_j}}\right)
$$
$$
 = \, {\psi}_*(\{\nabla_{\frac{\widetilde\partial}{\partial
\overline{z_j}}}\nabla_{\frac{\widetilde\partial}{\partial z_i}}
\widehat{t}_\alpha \, , \widehat{t}_\beta\}) \,+\,
(-1)^f\cdot
{\psi}_*(\{\nabla_{\frac{\widetilde\partial}{\partial z_i}}
\widehat{t}_\alpha \, ,
\nabla_{\frac{\widetilde\partial}{\partial z_j}}
\widehat{t}_\beta\}) \, ,  \leqno{(2.6)}
$$
where $\frac{\widetilde\partial}{\partial \overline{z_j}}$, as
before, denotes the lift of $\frac{\partial}{\partial
\overline{z_j}}$ using the chosen trivialization of the fibration,
and $f$ is the relative dimension. The holomorphicity of any
$\widehat{t}_\alpha$ yields the equality
$$
{\psi}_*(\{\nabla_{\frac{\widetilde\partial}{\partial \overline{z_j}}} 
\nabla_{\frac{\widetilde\partial}{\partial z_i}}
\widehat{t}_\alpha , \widehat{t}_\beta\}) \, =\, - \, {\psi}_*(C_h
(\frac{\widetilde\partial}{\partial z_i},
\frac{\widetilde\partial}{\partial \overline{z_j}})
\{\widehat{t}_\alpha \, , \widehat{t}_\beta\})\, .
$$

Now, for any $v =  \sum_{i=1}^m\sum_{\alpha = 1}^r v_{i,\alpha}
\frac{\partial}{\partial z_i}\otimes t_\alpha  \in
T_\xi X\otimes V_\xi$ we have
$$
\sum_{i=1}^m \sum_{j=1}^m \sum_{\alpha = 1}^r \sum_{\beta = 1}^r
v_{i,\alpha} \overline{v_{j,\beta}}\cdot {\psi}_*(C_h
(\frac{\widetilde\partial}{\partial z_i},
\frac{\widetilde\partial}{\partial \overline{z_j}})
\{\widehat{t}_\alpha , \widehat{t}_\beta\}) =
{\psi}_*(\sum _{i,j=1}^m C_h
(\frac{\widetilde\partial}{\partial z_i},
\frac{\widetilde\partial}{\partial \overline{z_j}})\{{\theta}_i,
\overline{{\theta}_j}\})\, ,
$$
where $\theta_i \,:=\, \sum_{\alpha =1}^r v_{i,a}\widehat{t}_a$
is the section of $K_{Y/X}\bigotimes L$ defined
on a neighborhood of $Y_\xi$. Let $e$ be a local section of $L$
defined around a point $y\in Y$ and $\{{\zeta}_1, {\zeta}_2,
\cdots , {\zeta}_f\}$ be a holomorphic
coordinate chart on $Y_\xi$ around $y$.
Set $\tau$ to be the local section of $TY$ defined around $y$
that satisfies the equality
$$
d{\zeta}_1\wedge {\zeta}_2\wedge
\cdots\wedge {\zeta}_f \otimes e \otimes \tau \, = \, \sum_{i=1}^m
\frac{\widetilde\partial}{\partial z_i}\otimes {\theta}_i
$$
of local sections of $K_{Y/X}\bigotimes L\bigotimes TY$. Using this
notation, the evaluation at $\xi$ of the function defined on a
neighborhood of $\xi$ by the fiber integral
$$
{\psi}_*(\{\nabla_{\frac{\widetilde\partial}{\partial
\overline{z_j}}}\nabla_{\frac{\widetilde\partial}{\partial z_i}}
\widehat{t}_\alpha \, , \widehat{t}_\beta\})
$$
in (2.6) coincides with
$$
{\psi}_*\left(\big\langle C_h(\tau ,\overline{\tau})e,e\big\rangle_h
\cdot d{\zeta}_1\wedge d{\zeta}_2\wedge 
\cdots \wedge d{\zeta}_f\wedge d\overline{{\zeta}_1}\wedge
d\overline{{\zeta}_2} \wedge \cdots \wedge
d\overline{{\zeta}_f}\right)\, . \leqno{(2.7)}
$$

The curvature $C_h$ is given to be positive. So using the
expression (2.7) we conclude that the evaluation at $\xi$
of the complex valued function around $\xi$ defined by
$$
\frac{1}{\sqrt{-1}}\cdot
\sum_{i=1}^m \sum_{j=1}^m \sum_{\alpha = 1}^r \sum_{\beta = 1}^r
v_{i,\alpha} \overline{v_{j,\beta}}\cdot
{\psi}_*(\{\nabla_{\frac{\widetilde\partial}{\partial \overline{z_j}}}
\nabla_{\frac{\widetilde\partial}{\partial z_i}}
\widehat{t}_\alpha \, , \widehat{t}_\beta\})
$$
is positive real; ${\psi}_*(\{\nabla_{\frac{\widetilde\partial}
{\partial \overline{z_j}}}
\nabla_{\frac{\widetilde\partial}{\partial z_i}}
\widehat{t}_\alpha \, , \widehat{t}_\beta\})$
is a term in (2.6). Consequently, to prove
the theorem it suffices to show that the last term in (2.6), namely
the complex valued function
$$
{\psi}_*(\{\nabla_{\frac{\widetilde\partial}{\partial z_i}}
\widehat{t}_\alpha \, ,
\nabla_{\frac{\widetilde\partial}{\partial z_j}}
\widehat{t}_\beta\}) \, ,  \leqno{(2.8)}
$$
which defined around $\xi$, actually vanishes at $\xi$.

Let $\omega$ denote the K\"ahler form on $Y_{\xi}$ obtained from
the curvature of the Chern connection $C_h$ on $L$.
The second part of the normal frame condition gives
$$
\big\langle (\nabla_{\frac{\widetilde\partial}{\partial z_i}}
\widehat{t}_\alpha\,)\big\vert_{Y_{\xi}},
\widehat{t}_\beta\big\vert_{Y_{\xi}}\big\rangle_{\omega}
\, = \, {\psi}_*(\{\nabla_{\frac{\widetilde\partial}{\partial z_i}}
\widehat{t}_\alpha , \widehat{t}_\beta\})(\xi) \, = \, 0
$$
for all $\alpha ,\beta \in [1,r]$ and $i\in [1,m]$;
here $\langle - , - \rangle_{\omega}$ is
the inner product defined using
$\omega$. We note that the above orthogonality condition does not
depend on the choice of the K\"ahler form on $Y_{\xi}$ that is
needed to define the orthogonality condition.

The above assertion that the smooth section
$$
\widetilde{t}_{i,\alpha}\, :=\,
(\nabla_{\frac{\widetilde\partial}{\partial z_i}}
\widehat{t}_\alpha)\big\vert_{Y_{\xi}}
$$
is orthogonal to $H^0(Y_{\xi}, \,  (K_{Y/X}\bigotimes
L)\big\vert_{Y_{\xi}})$ is equivalent to the condition that
$\widetilde{t}_{i,\alpha}$ is orthogonal to the space
$H^{f,0}(Y_{\xi},\, L\big\vert_{Y_{\xi}})$ of
$\Delta_{Y_{\xi}}''$-harmonic 
$(f,0)$-forms with values in $L$; here $f$, as before,
is the relative dimension. This Laplacian $\Delta_{Y_{\xi}}''$
corresponds to the Dolbeault operator on $L\big\vert_{Y_{\xi}}$
endowed with the Hermitian metric $h\big\vert_{Y_{\xi}}$.

The above orthogonality condition implies that
the equality
$$
\widetilde{t}_{i,\alpha} \,=\, \Delta_{Y_{\xi}}'' G_{Y_{\xi}}''
(\widetilde{t}_{i,\alpha}) \,=\,
(D_{Y_{\xi}}'')^* D_{Y_{\xi}}''
G_{Y_{\xi}}'' (\widetilde{t}_{i,\alpha})
$$
is valid; here $G_{Y_{\xi}}''$ is the Green operator
corresponding to $\Delta_{Y_{\xi}}''$. Therefore, we have the
following equality
$$
{\psi}_*(\{\nabla_{\frac{\widetilde\partial}{\partial z_i}}
\widehat{t}_\alpha , \nabla_{\frac{\widetilde\partial}{\partial z_j}}
\widehat{t}_\beta\})
 =  \int_{Y_{\xi}}
\langle D_{Y_{\xi}}'' (\widetilde{t}_{i,\alpha}),
 D_{Y_{\xi}}'' G_{Y_{\xi}}'' (\widetilde{t}_{j,\beta})\rangle_{\omega}
=  \int_{Y_{\xi}}
\langle (D_{Y_{\xi}}'')^* D_{Y_{\xi}}'' (\widetilde{t}_{i,\alpha}),
G_{Y_{\xi}}'' (\widetilde{t}_{j,\beta})\rangle_{\omega}
$$
for the function defined by (2.8).

Now, to prove that the function defined in (2.8) vanishes at
$\xi$ it suffices to establish the equality
$$
D_{Y_{\xi}}'' (\widetilde{t}_{i,\alpha}) \, = \, 0 \leqno{(2.9)}
$$
for all $i\in [1,m]$ and $\alpha \in [1,r]$.

To prove this we first note that
$$
D_{Y_{\xi}}'' (\widetilde{t}_{i,\alpha}) =
D_{Y_{\xi}}'' ((\nabla_{\frac{\widetilde
\partial}{\partial z_i}} \widehat{t}_\alpha\,)\big\vert_{Y_{\xi}}) =
(\sum_{k =1}^f
{\nabla}_{\frac{\partial}{\partial \overline{\zeta_k}}}
(\nabla_{\frac{\widetilde \partial}{\partial z_i}}
\widehat{t}_\alpha\,) d \overline{\zeta_k})\big\vert_{Y_{\xi}}
= -\sum_{k =1}^f C_h(\frac{\widetilde \partial}{\partial z_i},
\frac{\partial}{\partial \overline{\zeta_k}})
\widehat{t}_\alpha d \overline{\zeta_k}\, ,
$$
where $\{{\zeta}_1, \cdots , {\zeta}_{f-1}, {\zeta}_f\}$, as before
is a local holomorphic coordinate chart on $Y_{\xi}$.

The locally defined $(0,1)$-form
$\sum_{k =1}^f C_h(\frac{\widetilde \partial}{\partial z_i},
\frac{\partial}{\partial \overline{\zeta_k}})
d \overline{\zeta_k}$ is easily seen to be independent of the
coordinate function $\{{\zeta}_1, \cdots ,{\zeta}_f\}$. Indeed,
this locally defined $(0,1)$-form coincides with the contraction
$$
i_{\frac{\widetilde \partial}{\partial z_i}} C_h 
\vert_{Y_\xi}
$$
of the $(1,1)$-form $C_h$ with the vector field
$\frac{\widetilde \partial}{\partial z_i}$.

We have $H^{0,1}(Y_{\xi}) =0$, as, by
assumption, $H^1(Y_{\xi},\, {\Bbb Q}) =0$ and $Y_{\xi}$
is K\"ahler. Consequently, there is no
nonzero harmonic form of type $(0,1)$ on $Y_{\xi}$. Therefore, the
equality (2.9) is an immediate consequence of the following lemma.

\medskip
\noindent {\bf Lemma\, 2.10.}\, {\it The $(0,1)$-form
$i_{\frac{\widetilde \partial}{\partial z_i}} C_h 
\vert_{Y_\xi}$ on $Y_{\xi}$ is
harmonic.}
\medskip

{\it Proof of Lemma 2.10.}\, This form
$i_{\frac{\widetilde \partial}{\partial z_i}} C_h 
\vert_{Y_\xi}$ is evidently
$\overline{\partial}_{Y_{\xi}}$-closed.

To prove that it is also $\overline{\partial}^*_{Y_{\xi}}$-closed,
first note that the K\"ahler identity gives
$$
\sqrt{-1}\cdot
\overline{\partial}^*_{Y_{\xi}}\left(\sum_{k =1}^f C_h(\frac{
\widetilde\partial}{\partial z_i}\, ,
\frac{\partial}{\partial \overline{\zeta_k}})
d \overline{\zeta_k}\right) \,=\, {\Lambda}_{\omega}
{\partial}_{Y_{\xi}}\sum_{k =1}^f C_h(\frac{\widetilde
\partial}{\partial z_i}\, ,
\frac{\partial}{\partial \overline{\zeta_k}})
d \overline{\zeta_k}\, .
$$
The right-hand side coincides with
$$
{\Lambda}_{\omega}
\sum_{l =1}^f \sum_{k =1}^f \frac{\partial}{\partial{\zeta}_l}
C_h(\frac{\widetilde \partial}{\partial z_i}\, ,
\frac{\partial}{\partial\overline{\zeta_k}})
d {\zeta_l}\wedge d \overline{\zeta_k} \, =\,
{\Lambda}_{\omega}
\sum_{l=1}^f \sum_{k =1}^f\frac{\widetilde \partial}{\partial z_i}
C_h(\frac{\partial}{\partial{\zeta}_l}\, , \frac{\partial}{\partial
\overline{\zeta_k}})
d {\zeta_l}\wedge d \overline{\zeta_k}\, .
$$
The last equality is a obtained using the Bianchi identity. Since
${\Lambda}_{\omega}$ and $\frac{\widetilde \partial}{\partial z_i}$
commute, the following equality is obtained~:
$$
{\Lambda}_{\omega}
\sum_{l =1}^f \sum_{k =1}^f\frac{\widetilde \partial}{\partial z_i}   
C_h(\frac{\partial}{\partial{\zeta}_l}\, , \frac{\partial}{\partial
\overline{\zeta_k}})
d {\zeta_l}\wedge d \overline{\zeta_k} \, = \,
\frac{\widetilde \partial}{\partial z_i}({\Lambda}_{\omega} \omega)
\, =\, 0\, .
$$
This completes the proof of the lemma.$\hfill{\Box}$
\medskip

We already noted that the given condition $H^{1}(Y_{\xi},\, {\Bbb Q})
=0$ and the above lemma together imply the equality (2.9). This
completes the proof of the theorem.$\hfill{\Box}$
\medskip

If $H^{1}(Y_{\xi},\, {\Bbb Q})
\neq 0$, then the Picard group of $Y_{\xi}$ has continuous part.
If $\psi$ gives a locally trivial holomorphic fibration,
then the family of line bundle
$L\big\vert_{Y_x}$ gives an infinitesimal deformation map
$\tau : T_\xi X \, \longrightarrow \,
H^1(Y_{\xi}, {\cal O}_{Y_{\xi}})$.
It is easy to check that the harmonic $(0,1)$-form
$\sum_{k =1}^f C_h(\frac{\widetilde \partial}{\partial z_i},
\frac{\partial}{\partial \overline{\zeta_k}})
d \overline{\zeta_k}$ in Lemma 2.10 is the harmonic representative
of the image of the tangent vector
$\frac{\partial}{\partial z_i}$ under the homomorphism $\tau$.
(See \cite{ST}, \cite{BS} for a similar argument.)

\section{Applications of the positivity of direct images}

Let $E$ be an ample vector bundle of rank $r$
on a projective manifold $X$. Take
$\lambda \,= \,({\lambda_1}, {\lambda_2}, \cdots , {\lambda_r}) 
\,\in\, {\Bbb N}^r$
such that ${\lambda_i} \geq {\lambda_j}$ if $i\leq j$. Let
${\Gamma}^{\lambda} E$ denote the vector bundle associated to $E$ for
the weight $\lambda$. If, in particular, ${\lambda}_i =0$ for
$i\geq 2$, then ${\Gamma}^{\lambda} E$ is the symmetric power
$S^{{\lambda}_1}(E)$; if ${\lambda}_i =0$ for $i\geq k+1$ and
${\lambda}_i =1$ for $i\leq k$, then ${\Gamma}^{\lambda} E$ is the
exterior power ${\bigwedge}^{k}(E)$. Let
$$
\psi \, :\, M_{\lambda}(E) \,\longrightarrow\, X
$$
denote the associated flag bundle over $X$ and
$L_{\lambda}\longrightarrow M_{\lambda}(E)$ is the
corresponding line bundle. If ${\lambda}_2 =0$, then $M_{\lambda}(E) =
{\Bbb P}(E)$ and $L_{\lambda}= {\cal O}_{{\Bbb P}(E)}({\lambda_1}+r)$.
The direct image ${\psi}_* (L_{\lambda}\bigotimes
K_{M_{\lambda}(E)/X})$ coincides with ${\Gamma}^{\lambda} E\bigotimes
({\bigwedge}^r E)^{\otimes l}$,
where $l \in [1,r]$ such that ${\lambda}_l
\neq 0$ and ${\lambda}_{l+1} =0$;
$K_{M_{\lambda}(E)/X}$ is the relative
canonical bundle for the projection $\psi$. In this special setup
Theorem 2.3 reads as follows~:

\medskip
\noindent {\bf Theorem\, 3.1.}\, {\it The vector bundle
${\Gamma}^{\lambda}E\bigotimes ({\bigwedge}^r E)^{\otimes l}$
is Nakano positive. In particular,
setting $\lambda = (k,0,\cdots ,0,0)$ the vector bundle $S^k(E)
\bigotimes\det E$ is Nakano positive.}
\medskip

Theorem 3.1 combined with a vanishing theorem of
Nakano proved in \cite{Na}
immediately gives as a corollary the following result of Demailly
proved in \cite{De}.

\medskip
\noindent {\bf Corollary\, 3.2.}\, {\it Let $E$ be an ample vector
bundle on a projective manifold $X$ of dimension $n$. Then
$$
H^{n,i}(X,\, {\Gamma}^{\lambda} E\otimes (\det E)^{\otimes l})\, =\, 0 
$$
for $i\geq 1$ and $\lambda$ as above.}
\medskip

We note that the special case Corollary 3.2 where
${\Gamma}^{\lambda} E
= S^k(E)$ was proved by Griffiths \cite{Gr}.

\medskip
\noindent {\bf Remark\, 3.3.}\, The curvature terms of
the $L^2$-metric on the vector bundle
${\Gamma}^{\lambda} E\otimes (\det E)^{\otimes l}$ are
$$
{\psi}_*(\{\nabla_{\frac{\widetilde\partial}{\partial
\overline{z_j}}}\nabla_{\frac{\widetilde\partial}{\partial z_i}}
\widehat{t}_\alpha \, , \widehat{t}_\beta\})
$$
The curvature of the dual metric is the
negative of the transpose of the initial curvature.
Hence, the sesquilinear form on $TX\otimes E^\star$ computed with
the dual metric is Nakano negative. As an immediate consequence,
for an ample vector bundle $E$ we have
$$
H^{p,n}(X,\, {\Gamma}^{\lambda} E\otimes (\det E)^{\otimes l})
\, =\, 0
$$
if $p\geq 1$ and $\lambda$ as above.
\medskip

Let $X$ be a compact complex manifold equipped with
a Hermitian structure $\omega$. In \cite{DPS}, the
notion of a {\it numerically effective line bundle} on
$X$ is defined to be a holomorphic line bundle $L$ satisfying
the condition that given any $\epsilon > 0$, there is a Hermitian
metric $h_{\epsilon}$ on $L$ such that
$$
{\Theta}_{h_{\epsilon}} \, \geq  \, - \,\epsilon \cdot \omega \, ,
$$
where ${\Theta}_{h_{\epsilon}}$ is the Chern curvature for
the Hermitian metric $h_{\epsilon}$ \cite[Definition 1.2]{DPS}.
The manifold $X$ being compact, this condition,
of course, does not depend on $\omega$. A
vector bundle $E$ over $X$ is called {\it numerically effective}
if the tautological line bundle ${\cal O}_{{\Bbb P}(E)}(1)$ over
${\Bbb P}(E)$ is numerically effective \cite[Definition 1.9]{DPS}.

We have the following proposition as an easy consequence of the
proof of Theorem 2.3.

\medskip
\noindent {\bf Proposition\, 3.4.}\, {\it Let $X$ be a compact
complex manifold equipped with a Hermitian form $\omega$.
A vector bundle $E$ over
$X$ is numerically effective if and only if there is a Hermitian
metric $h_{k,\varepsilon}$ on $S^k(E)\det E$, for all $k\geq 1$
and all $\varepsilon > 0$, such that
$$
\sqrt{-1}\cdot C_{h_{k,\varepsilon}}(S^k(E)\otimes \det E)\, >
\, - \varepsilon\cdot\omega\otimes Id_{S^k(E)\otimes \det E} \, ,
$$
where $C_{h_{k,\varepsilon}}(S^k(E)\otimes \det E)$
denotes the curvature of the metric $h_{k,\varepsilon}$.
The inequality is in the sense of Nakano.}
\medskip

{\it Proof.}\, Let $E$ be a numerically effective vector
bundle of rank $r$
over the compact complex manifold $X$. Let $Y$ denote the
projective bundle ${\Bbb P}(E)$ over $X$. The natural projection
of $Y$ to $X$ will be denoted by $\psi$.

The tautological line bundle ${\cal O}_{{\Bbb P}(E)}(1)$ on
${\Bbb P}(E)$ is numerically effective.
For any $\varepsilon > 0$, consider the Hermitian metric
$h_{\varepsilon}$ on ${\cal O}_{{\Bbb P}(E)}(1)$ as in
Definition 1.2 (pp. 299) of \cite{DPS}. This Hermitian metric
$h_{\varepsilon}$ induces a Hermitian metric on each
${\cal O}_{{\Bbb P}(E)}(k)$, where $k > 1$.

Consequently,
we have a Hermitian metric $h_{k,\varepsilon}$ on each $S^k(E)
\bigotimes\det E$
obtained as the $L^2$-metric for the above
Hermitian metric on ${\cal O}_{{\Bbb P}(E)}(k+r)$ obtained from
$h_{\varepsilon}$. Now, from the
proof of Theorem 2.3 it can be deduced that the
condition for numerically effectiveness of a line bundle
in Definition 1.2 (pp. 299) of \cite{DPS},
given in terms of Hermitian metrics $h_{\varepsilon}$,
ensures that the inequality in the proposition is valid.

Conversely, let $E$ be a vector bundle such that the inequality
condition in the statement of the proposition is valid. Now it
follows immediately from the criterion of numerically effectiveness
given in Theorem 1.12 (pp. 306) of \cite{DPS} that such a vector
bundle $E$ must be numerically effective. This completes the
proof of the proposition.$\hfill{\Box}$
\medskip

\noindent {\bf Remark\, 3.5.}\, Let $E$ be a vector bundle of rank $r$
and set $L\, =\, {\cal O}_{{\Bbb P}(E)}(r+k)$ in (2.6). The
expression for the curvature
of the $L^2$-metric on $S^k E\otimes \det E$ shows that
if $E$ is ample, then for each $k\,\geq \,1$, there is a
Hermitian metrics on $S^k E\otimes \det E$ such that
there exits a positive number $\varepsilon \, >\, 0$ with the
property that for every positive integer $k$, the inequality
$$
\sqrt{-1}\cdot c(S^k E\otimes \det E) \, > \, (k+r) \varepsilon
\omega\otimes {\rm Id}_{S^k E\otimes \det E}
$$
is valid, where $\omega$ is any fixed Hermitian form on $X$; the
inequality is in the sense of Nakano.

We observe that the above property is actually a characterization
of ampleness. Indeed, if $E$ is a vector bundle which satisfies the
above condition, then fix any metric on the line bundle $\det E$.
By subtracting it to $S^k E\otimes \det E$, 
we get a metric on $S^k E$
whose curvature is, for a large $k$, Nakano positive.
Hence, $S^k E$ must be ample. This immediately yields the
ampleness of $E$.



\begin{thebibliography}{AAAA}

\bibitem[BS]{BS} I. Biswas and G. Schumacher : Determinant bundle,
Quillen metric, and Petersson-Weil form on moduli spaces.
Geom. Funct. Anal. {\bf 9} (1999), 226--255.

\bibitem[De]{De} J.-P. Demailly : Vanishing theorems for tensor
powers of an ample vector bundle. Invent. Math. {\bf 91} (1988),
203--220.

\bibitem[DS]{DS} J.-P. Demailly and H. Skoda : Relations entre les
notions de positivit\'e de P.A. Griffiths et de S. Nakano.
S\'eminaire Pierre Lelong -- Henri Skoda (Analyse) Ann\'ees 1978/79.
304--309, Lecture Notes in Math., 822, Springer, Berlin, 1980.

\bibitem[DPS]{DPS} J.-P. Demailly, T. Peternell and M. Schneider :
Compact complex manifolds with numerically effective tangent bundles.
Jour. Alg. Geom. {\bf 3} (1994), 295--345.

\bibitem[Gr]{Gr} P.A. Griffiths : Hermitian differential geometry,
Chern classes,  and positive vector bundles. {\it Global Analysis},
Papers in Honor of K. Kodaira (Ed. D. C. Spencer and S. Iyanaga),
185--251. Princeton Univ. Press, 1969.

\bibitem[Ha]{Ha} R. Hartshorne : Ample vector bundles.
Inst. Hautes \'Etudes Sci. Publ. Math {\bf 29} (1966), 63--94.

\bibitem[Na]{Na} S. Nakano : On complex analytic vector bundles.
Jour. Math. Soc. Japan {\bf 7} (1955), 1--12.

\bibitem[Mo]{Mo} Ch. Mourougane : Images directes de fibr\'es en
droites adjoints. Publ. RIMS, Kyoto Univ. {\bf 33} (1997), 893--916.

\bibitem[ST]{ST} G. Schumacher and M. Toma : On the Petersson-Weil
metric for the moduli space of Hermite-Einstein bundles and its
curvature. Math. Ann. {\bf 293} (1992), 101--107.

\bibitem[Um]{Um} H. Umemura : Some results in the theory of
vector bundles. Nagoya Math. Jour. {\bf 52} (1973), 97--128.

\end{thebibliography}
\end{document}